\documentclass[11pt]{amsart}

\theoremstyle{plain}
\newtheorem{theorem}                {Theorem}      [section]
\newtheorem{proposition}  [theorem]  {Proposition}

\newtheorem{lemma}        [theorem]  {Lemma}

\theoremstyle{definition}

\newtheorem{remark}       [theorem]  {Remark}
\newtheorem{definition}   [theorem]  {Definition}

\def \s{\mbox{${\mathbb S}$}}

\DeclareMathOperator{\trace}{trace}

\DeclareMathOperator{\Span}{span}
 
\DeclareMathOperator{\cst}{constant}

\numberwithin{equation}{section}

\begin{document}

\title[Biharmonic submanifolds in Sasakian space forms]{Explicit formulas for biharmonic submanifolds in Sasakian space forms}

\author{D.~Fetcu}
\author{C.~Oniciuc}

\address{Department of Mathematics\\
"Gh. Asachi" Technical University of Iasi\\
Bd. Carol I no. 11 \\
700506 Iasi, Romania} \email{dfetcu@math.tuiasi.ro}

\address{Faculty of Mathematics\\ "Al.I. Cuza" University of Iasi\\
Bd. Carol I no. 11 \\
700506 Iasi, Romania} \email{oniciucc@uaic.ro}

\begin{abstract} We classify the biharmonic Legendre curves in a
Sasakian space form, and obtain their explicit parametric equations
in the $(2n+1)$-dimensional unit sphere endowed with the canonical
and deformed Sasakian structures defined by Tanno. Then, composing
with the flow of the Reeb vector field, we transform a biharmonic
integral submanifold into a biharmonic anti-invariant submanifold.
Using this method we obtain new examples of biharmonic submanifolds
in spheres and, in particular, in $\mathbb{S}^{7}$.

\end{abstract}

\date{}



\dedicatory{Dedicated to Professor Neculai Papaghiuc on his 60-th
birthday}

\thanks{The authors were partially supported by the Grant CEEX, ET, 5871/2006, Romania.}

\maketitle

\section{Introduction}
\setcounter{equation}{0}

\textit{Biharmonic maps} between Riemannian manifolds
$\phi:(M,g)\to(N,h)$  are the critical points of the
\textit{bienergy functional}
$E_{2}(\phi)=\frac{1}{2}\int_{M}|\tau(\phi)|^{2} \ v_{g}$. They
represent a natural generalization of the well-known
\textit{harmonic maps} (\cite{Eells}), the critical points of the
\textit{energy functional} $E(\phi)=\frac{1}{2}\int_{M}|d\phi|^{2} \
v_{g}$, and of the biharmonic submanifolds in Euclidean spaces
defined by B.-Y. Chen (\cite{Chen}).

The Euler-Lagrange equation for the energy functional is
$\tau(\phi)=0$, where $\tau(\phi)=\trace\nabla d\phi$ is the
tension field, and the Euler-Lagrange equation for the bienergy
functional was derived by G. Y. Jiang in \cite{Jiang}:
$$
\begin{array}{cl}
\tau_{2}(\phi)&=-\Delta\tau(\phi)-\trace\
R^{N}(d\phi,\tau(\phi))d\phi\\ \\
&=0.\end{array}
$$
Since any harmonic map is biharmonic, we are interested in
non-harmonic biharmonic maps, which are called
\textit{proper-biharmonic}.

A special case of biharmonic maps is represented by the biharmonic
Riemannian immersions, or biharmonic submanifolds. There are several
results of classification or construction for such submanifolds in
space forms (\cite{MontaldoOniciuc}, \cite{BMO}). Then, the next
step would be the study of biharmonic submanifolds in Sasakian space
forms. In this context J. Inoguchi classified in \cite{Ino} the
proper-biharmonic Legendre curves and Hopf cylinders in a
3-dimensional Sasakian space form $M^{3}(c)$, and in
\cite{FetcuOniciuc} the explicit parametric equations were obtained.
Then, T. Sasahara and his collaborators studied the biharmonic
integral surfaces and 3-dimensional biharmonic anti-invariant
submanifolds in $\mathbb{S}^{5}$ (\cite{Sasahara2},
\cite{Sasahara1}).

Recent results on biharmonic submanifolds in spaces of nonconstant
sectional curvature were obtained by T. Ichiyama, J. Inoguchi and H.
Urakawa in \cite{Ura}, by Y.-L. Ou and Z.-P. Wang in \cite{Ou}, and
by W. Zhang in \cite{Zhang}.

Biharmonic submanifolds in pseudo-Euclidean spaces were also
studied, and many examples and classification results were obtained
(for example, see \cite{Arv}, \cite{Chen}).

The goals of our paper are to obtain new classification results for
biharmonic Legendre curves in any dimensional Sasakian space form
and to provide a method for constructing biharmonic submanifolds. In
order to obtain explicit examples, we use the $(2n+1)$-dimensional
unit sphere $\mathbb{S}^{2n+1}$ as a model of Sasakian space form.

For a general account of biharmonic maps see \cite{MontaldoOniciuc}
and \textit{The Bibliography of Biharmonic Maps} \cite{bibl}.

\noindent \textbf{Conventions.} We work in the $C^{\infty}$
category, that means manifolds, metrics, connections and maps are
smooth. The Lie algebra of the vector fields on $M$ is denoted by
$C(TM)$.

\section{Preliminaries}
\setcounter{equation}{0}

\subsection{Contact manifolds}

A \textit{contact metric structure} on a manifold $N^{2n+1}$ is
given by $(\varphi,\xi,\eta,g)$, where $\varphi$ is a tensor field
of type $(1,1)$ on $N$, $\xi$ is a vector field, $\eta$ is an 1-form
and $g$ is a Riemannian metric such that
$$
\begin{array}{c} \varphi^{2}=-I+\eta\otimes\xi,\ \
\eta(\xi)=1,\\ \\ g(\varphi X,\varphi Y)=g(X,Y)-\eta(X)\eta(Y),\ \
\ g(X,\varphi Y)=d\eta(X,Y),\ \ \forall X,Y\in C(TN).\end{array}
$$
A contact metric structure $(\varphi,\xi,\eta,g)$ is called {\it
normal} if
$$
N_{\varphi}+2d\eta\otimes\xi=0,
$$
where
$$
N_{\varphi}(X,Y)=[\varphi X,\varphi Y]-\varphi \lbrack \varphi
X,Y]-\varphi \lbrack X,\varphi Y]+\varphi ^{2}[X,Y],\ \ \forall
X,Y\in C(TN),
$$
is the Nijenhuis tensor field of $\varphi$.

\noindent A contact metric manifold $(N,\varphi,\xi,\eta,g)$ is a
\textit{Sasakian manifold} if it is normal or, equivalently, if
$$
(\nabla_{X}\varphi)(Y)=g(X,Y)\xi-\eta(Y)X,\ \ \forall X,Y\in
C(TN).
$$

\noindent The {\it contact distribution} of a Sasakian manifold
$(N,\varphi,\xi,\eta,g)$ is defined by $\{X\in TN:\eta(X)=0\}$, and
an integral curve of the contact distribution is called {\it
Legendre curve}. A submanifold $M$ of $N$ which is tangent to $\xi$
is said to be {\it anti-invariant} if $\varphi$ maps any vector
tangent to $M$ and normal to $\xi$ to a vector normal to $M$.

\noindent Let $(N,\varphi,\xi,\eta,g)$ be a Sasakian manifold. The
sectional curvature of a 2-plane generated by $X$ and $\varphi X$,
where $X$ is an unit vector orthogonal to $\xi$, is called
\textit{$\varphi$-sectional curvature} determined by $X$. A
Sasakian manifold with constant $\varphi$-sectional curvature $c$
is called a \textit{Sasakian space form} and it is denoted by
$N(c)$.

\noindent The curvature tensor field of a Sasakian space form
$N(c)$ is given by
\begin{equation}\label{eccurv}
\begin{array}{ll}
R(X,Y)Z=&\frac{c+3}{4}\{g(Z,Y)X-g(Z,X)Y\}+\frac{c-1}{4}\{\eta(Z)\eta(X)Y-\\ \\
&-\eta(Z)\eta(Y)X+g(Z,X)\eta(Y)\xi-g(Z,Y)\eta(X)\xi+\\
\\&+g(Z,\varphi Y)\varphi X-g(Z,\varphi X)\varphi
Y+2g(X,\varphi Y)\varphi Z\}.
\end{array}
\end{equation}

Let $\mathbb{S}^{2n+1}=\{z\in\mathbb{C}^{n+1}: \vert z\vert=1\}$
be the unit $(2n+1)$-dimensional sphere endowed with its standard
metric field $g_{0}$. Consider the following structure tensor
fields on $\mathbb{S}^{2n+1}$: $\xi_{0}=-\mathcal{I}z$ for each
$z\in \mathbb{S}^{2n+1}$, where $\mathcal{I}$ is the usual almost
complex structure on $\mathbb{C}^{n+1}$ defined by
$$
\mathcal{I}z=(-y^{1},...,-y^{n+1},x^{1},...,x^{n+1}),
$$
for $z=(x^{1},...,x^{n+1},y^{1},...,y^{n+1})$, and
$\varphi_{0}=s\circ \mathcal{I}$, where $s:T_{z}\mathbb{C}^{n+1}\to
T_{z}\mathbb{S}^{2n+1}$ denotes the orthogonal projection. Equipped
with these tensors, $\mathbb{S}^{2n+1}$ becomes a Sasakian space
form with $\varphi_{0}$-sectional curvature equal to 1 (\cite{BB}).

\noindent Now, consider the deformed structure on
$\mathbb{S}^{2n+1}$, introduced by Tanno in [13],
$$
\eta=a\eta_{0}, \ \ \xi=\frac{1}{a}\xi_{0},\ \
\varphi=\varphi_{0},\ \ g=a g_{0}+a(a-1)\eta_{0}\otimes\eta_{0},
$$
where $a$ is a positive constant. The structure
$(\varphi,\xi,\eta,g)$ is still a Sasakian structure and
$(\mathbb{S}^{2n+1},\varphi,\xi,\eta,g)$ is a Sasakian space form
with constant $\varphi$-sectional curvature $c=\frac{4}{a}-3$
([1]).

\subsection{3-Sasakian manifolds}

If a manifold $N$ admits three Sasakian structures
$(\varphi_{a},\xi_{a},\eta_{a},g)$, $a=1,2,3$, satisfying
$$
\begin{array}{c}
\varphi_{c}=-\varphi_{a}\varphi_{b}+\eta_{b}\otimes\xi_{a}=\varphi_{b}\varphi_{a}-\eta_{a}\otimes\xi_{b},\\
\\
\xi_{c}=-\varphi_{a}\xi_{b}=\varphi_{b}\xi_{a},\ \ \
\eta_{c}=-\eta_{a}\circ\varphi_{b}=\eta_{b}\circ\varphi_{a},
\end{array}
$$
for an even permutation $(a,b,c)$ of $(1,2,3)$, then the manifold is
said to have a \textit{Sasakian 3-structure}. The dimension of such
a manifold is of the form $4n+3$. We note that the maximum dimension
of a submanifold of a 3-Sasakian manifold $N^{4n+3}$ which is an
integral submanifold with respect to all three Sasakian structures
is $n$.

\section{Biharmonic Legendre curves in Sasakian space forms}
\setcounter{equation}{0}

\begin{definition}
Let $(N^{m},g)$ be a Riemannian manifold and $\gamma:I\to N$ a curve
parametrized by arc length, that is $\vert\gamma'\vert=1$. Then
$\gamma$ is called a \textit{Frenet curve of osculating order r},
$1\leq r\leq m$, if there exists orthonormal vector fields
$E_{1},E_{2},...,E_{r}$ along $\gamma$ such that
$$
\left\{\begin{array}{lc}E_{1}=\gamma'=T\\ \\
\nabla_{T}E_{1}=\kappa_{1}E_{2}\\
\\\nabla_{T}E_{2}=-\kappa_{1}E_{1}+\kappa_{2}E_{3}\\...\\ \\\nabla_{T}E_{r}=-\kappa_{r-1}E_{r-1}\end{array}\right.,
$$
where $\kappa_{1},...,\kappa_{r-1}$ are positive functions on $I$.
\end{definition}

\begin{remark}A geodesic is a Frenet curve of osculating order 1;
a \textit{circle} is a Frenet curve of osculating order 2 with
$\kappa_{1}=\cst$; a \textit{helix of order r}, $r\geq 3$, is a
Frenet curve of osculating order $r$ with
$\kappa_{1},...,\kappa_{r-1}$ constants; a helix of order $3$ is
called, simply, helix.
\end{remark}

Now let $(N^{2n+1},\varphi,\xi,\eta,g)$ be a Sasakian space form
with constant $\varphi$-sectional curvature $c$ and $\gamma:I\to N$
a Legendre Frenet curve of osculating order $r$. As
\begin{equation}
\begin{array}{llc}
\nabla_{T}^{3}T=&(-3\kappa_{1}\kappa_{1}')E_{1}+(\kappa_{1}''-\kappa_{1}^{3}-\kappa_{1}\kappa_{2}^{2})E_{2}+(2\kappa_{1}'\kappa_{2}+\kappa_{1}\kappa_{2}')E_{3}\\
\\
&+\kappa_{1}\kappa_{2}\kappa_{3}E_{4}
\end{array}
\end{equation}
and
\begin{equation}
R(T,\nabla_{T}T)T=-\frac{(c+3)\kappa_{1}}{4}E_{2}-\frac{3(c-1)\kappa_{1}}{4}g(E_{2},\varphi
T)\varphi T,
\end{equation}
we get
\begin{equation}\label{ec-bih-explicit}
\begin{array}{rl}
\tau_{2}(\gamma)=&\nabla_{T}^{3}T-R(T,\nabla_{T}T)T\\ \\
=&(-3\kappa_{1}\kappa_{1}')E_{1}+\Big(\kappa_{1}''-\kappa_{1}^{3}-\kappa_{1}\kappa_{2}^{2}+\frac{(c+3)\kappa_{1}}{4}\Big)E_{2}\\
\\&+(2\kappa_{1}'\kappa_{2}+\kappa_{1}\kappa_{2}')E_{3}+\kappa_{1}\kappa_{2}\kappa_{3}E_{4}+\frac{3(c-1)\kappa_{1}}{4}g(E_{2},\varphi T)\varphi T.
\end{array}
\end{equation}

\noindent In the following we shall solve the biharmonic equation
$\tau_{2}(\gamma)=0$. The problem is to find the relation between
$\varphi T$ and the Frenet frame field. The simplest two cases are
provided by $\frac{3(c-1)\kappa_{1}}{4}g(E_{2},\varphi T)=0$. So,

\noindent \textbf{Case I: $\mathbf{c=1}$.}

\noindent In this case $\gamma$ is proper-biharmonic if and only
if
$$
\left\{\begin{array}{l} \kappa_{1}=\cst>0,\ \ \kappa_{2}=\cst\\ \\
\kappa_{1}^{2}+\kappa_{2}^{2}=1\\ \\
\kappa_{2}\kappa_{3}=0\end{array}\right..
$$

\noindent One obtains
\begin{theorem}\label{teocase1}If $c=1$ and $n\geq 2$, then $\gamma$ is
proper-biharmonic if and only if either $\gamma$ is a circle with
$\kappa_{1}=1$, or $\gamma$ is a helix with
$\kappa_{1}^{2}+\kappa_{2}^{2}=1$.
\end{theorem}

\begin{remark} If $n=1$ and $\gamma$ is a non-geodesic Legendre curve we have $\nabla_{T}T=\pm\kappa_{1}\varphi
T$ and then $E_{2}=\pm\varphi T$ and
$\nabla_{T}E_{2}=\pm\nabla_{T}\varphi
T=\pm(\xi\mp\kappa_{1}T)=-\kappa_{1}T\pm\xi$. Therefore
$\kappa_{2}=1$ and $\gamma$ cannot be biharmonic.
\end{remark}

\noindent\textbf{Case II: $\mathbf{c\neq 1,\ E_{2}\perp\varphi
T}$.}

\noindent In this case $\gamma$ is proper-biharmonic if and only
if
$$
\left\{\begin{array}{l} \kappa_{1}=\cst>0,\ \ \kappa_{2}=\cst\\ \\
\kappa_{1}^{2}+\kappa_{2}^{2}=\frac{c+3}{4}\\ \\
\kappa_{2}\kappa_{3}=0\end{array}\right..
$$

\noindent Before stating the theorem we need the following

\begin{lemma} Let $\gamma$ be a Legendre Frenet curve of
osculating order 3 and $E_{2}\perp\varphi T$. Then
$\{T=E_{1},E_{2},E_{3},\varphi T,\xi,\nabla_{T}\varphi T\}$ are
linearly independent, in any point, and hence $n\geq 3$.
\end{lemma}

\begin{proof} Since $\gamma$ is a Frenet curve of osculating order
3, we have
$$
\left\{\begin{array}{lc}E_{1}=\gamma'=T\\ \\
\nabla_{T}E_{1}=\kappa_{1}E_{2}\\
\\\nabla_{T}E_{2}=-\kappa_{1}E_{1}+\kappa_{2}E_{3}\\
\\\nabla_{T}E_{3}=-\kappa_{2}E_{2}\end{array}\right..
$$

\noindent It is easy to see that, in an arbitrary point, the
system
$$
S_{1}=\{T=E_{1},E_{2},E_{3},\varphi T,\xi,\nabla_{T}\varphi T\}
$$
has only non-zero vectors and
$$
T\perp E_{2},\ \ T\perp E_{3},\ \ T\perp\varphi T,\ \ T\perp\xi,\
\ T\perp\nabla_{T}\varphi T.
$$
Thus $S_{1}$ is linearly independent if and only if $
S_{2}=\{E_{2},E_{3},\varphi T,\xi,\nabla_{T}\varphi T\}$ is linearly
independent. Further, as
$$
E_{2}\perp\xi,\ E_{2}\perp\nabla_{T}\varphi T,\ E_{3}\perp\xi,\
E_{3}\perp\nabla_{T}\varphi T,\ \varphi T\perp\xi,\ \varphi
T\perp\nabla_{T}\varphi T,
$$
and
$$
\ E_{2}\perp E_{3}\perp \varphi T,
$$
it follows that $S_{2}$ is linearly independent if and only if $
S_{3}=\{\xi,\nabla_{T}\varphi T\}$ is linearly independent. But
$\nabla_{T}\varphi T=\xi+\kappa_{1}\varphi E_{2}$, $\kappa_{1}\neq
0$, and therefore $S_{3}$ is linearly independent.
\end{proof}

\noindent Now we can state

\begin{theorem}Assume that $c\neq 1$ and $\nabla_{T}T\perp\varphi
T$. We have

\noindent 1) If $c\leq -3$ then $\gamma$ is biharmonic if and only
if it is a geodesic.

\noindent 2) If $c>-3$ then $\gamma$ is proper-biharmonic if and
only if either

a) $n\geq 2$ and $\gamma$ is a circle with
$\kappa_{1}^{2}=\frac{c+3}{4}$. In this case $\{E_{1},E_{2},\varphi
T,\xi\}$ are linearly independent,

or

b) $n\geq 3$ and $\gamma$ is a helix with
$\kappa_{1}^{2}+\kappa_{2}^{2}=\frac{c+3}{4}$. In this case
$\{E_{1},E_{2},E_{3},\varphi T,$ $\xi,\nabla_{T}\varphi T\}$ are
linearly independent.

\end{theorem}

\noindent\textbf{Case III: $\mathbf{c\neq 1,\
E_{2}\parallel\varphi T}$.}

\noindent In this case $\gamma$ is proper-biharmonic if and only
if
$$
\left\{\begin{array}{l} \kappa_{1}=\cst>0,\ \ \kappa_{2}=\cst\\ \\
\kappa_{1}^{2}+\kappa_{2}^{2}=c\\ \\
\kappa_{2}\kappa_{3}=0\end{array}\right..
$$

\noindent We can assume that $E_{2}=\varphi T$. Then we have
$\nabla_{T}T=\kappa_{1}E_{2}=\kappa_{1}\varphi T$,
$\nabla_{T}E_{2}=\nabla_{T}\varphi T=\xi-\kappa_{1}T$. That means
$E_{3}=\xi$ and $\kappa_{2}=1$. Hence
$\nabla_{T}E_{3}=\nabla_{T}\xi=-\varphi T=-E_{2}$.

\noindent Therefore

\begin{theorem}If $c\neq 1$ and $\nabla_{T}T\parallel\varphi T$,
then $\{T,\varphi T,\xi\}$ is the Frenet frame field of $\gamma$
and we have

1) If $c\leq 1$ then $\gamma$ is biharmonic if and only if it is a
geodesic.

2) If $c>1$ then $\gamma$ is proper-biharmonic if and only if it
is a helix with $\kappa_{1}^{2}=c-1$ (and $\kappa_{2}=1$).
\end{theorem}

\begin{remark} If $n=1$ then $\nabla_{T}T\parallel \varphi T$ and
we reobtain Inoguchi's result \cite{Ino}.
\end{remark}

\noindent\textbf{Case IV: $\mathbf{c\neq 1}$ and
$\mathbf{g(E_{2},\varphi T)}$ is not constant $\mathbf{0,1}$ or
$\mathbf{-1}$.}

\noindent Assume that $\gamma$ is a Legendre Frenet curve of
osculating order $r$, $4\leq r$ $\leq 2n+1$, $n\geq 2$. If $\gamma$
is biharmonic it follows that $\varphi
T\in\Span\{E_{2},E_{3},E_{4}\}$.

\noindent Now, we denote $f(t)=g(E_{2},\varphi T)$ and
differentiating along $\gamma$ we obtain
$$
\begin{array}{ll}
f'(t)&=g(\nabla_{T}E_{2},\varphi T)+g(E_{2},\nabla_{T}\varphi
T)=g(\nabla_{T}E_{2},\varphi T)+g(E_{2},\xi+\kappa_{1}\varphi E_{2})\\
\\&=g(\nabla_{T}E_{2},\varphi T)=g(-\kappa_{1}T+\kappa_{2}E_{3},\varphi
T)\\ \\&=\kappa_{2}g(E_{3},\varphi T).
\end{array}
$$

\noindent Since $\varphi T=g(\varphi T,E_{2})E_{2}+g(\varphi
T,E_{3})E_{3}+g(\varphi T,E_{4})E_{4}$, the curve $\gamma$ is
proper-biharmonic if and only if
$$
\left\{\begin{array}{l}\kappa_{1}=\cst>0\\ \\
\kappa_{1}^{2}+\kappa_{2}^{2}=\frac{c+3}{4}+\frac{3(c-1)}{4}f^{2}\\
\\\kappa_{2}'=-\frac{3(c-1)}{4}fg(\varphi T,E_{3})\\
\\\kappa_{2}\kappa_{3}=-\frac{3(c-1)}{4}fg(\varphi T,E_{4})
\end{array}\right..
$$

\noindent Using the expression of $f'(t)$ we see that the third
equation of the above system is equivalent to
$$
\kappa_{2}^{2}=-\frac{3(c-1)}{4}f^{2}+\omega_{0},
$$
where $\omega_{0}=\cst$. Replacing in the second equation it
follows
$$
\kappa_{1}^{2}=\frac{c+3}{4}-\omega_{0}+\frac{3(c-1)}{2}f^{2},
$$
which implies $f=\cst$. Thus $\kappa_{2}=\cst>0$, $g(E_{3},\varphi
T)=0$ and then $\varphi T=fE_{2}+g(\varphi T,E_{4})E_{4}$. It
follows that there exists an unique constant $\alpha_{0}\in
(0,2\pi)\setminus\{\frac{\pi}{2},\pi,\frac{3\pi}{2}\}$ such that
$f=\cos\alpha_{0}$ and $g(\varphi T,E_{4})=\sin\alpha_{0}$.

\noindent We can state

\begin{theorem}Let $c\neq 1$, $n\geq 2$ and $\gamma$ a Legendre
Frenet curve of osculating order $r\geq 4$ such that
$g(E_{2},\varphi T)$ is not constant $0,1$ or $-1$. We have

a) If $c\leq -3$ then $\gamma$ is biharmonic if and only if it is
a geodesic.

b) If $c>-3$ then $\gamma$ is proper-biharmonic if and only if
$\varphi T=\cos\alpha_{0}E_{2}+\sin\alpha_{0}E_{4}$ and
$$
\left\{\begin{array}{l}\kappa_{1}=\cst>0,\ \ \kappa_{2}=\cst\\ \\
\kappa_{1}^{2}+\kappa_{2}^{2}=\frac{c+3}{4}+\frac{3(c-1)}{4}\cos^{2}\alpha_{0}\\
\\\kappa_{2}\kappa_{3}=-\frac{3(c-1)}{8}\sin 2\alpha_{0}
\end{array}\right.,
$$
where $\alpha_{0}\in
(0,2\pi)\setminus\{\frac{\pi}{2},\pi,\frac{3\pi}{2}\}$ is a
constant such that $c+3+3(c-1)\cos^{2}\alpha_{0}>0$ and
$3(c-1)\sin 2\alpha_{0}<0$.
\end{theorem}

\begin{remark}
In this case we may obtain biharmonic curves which are not
helices.
\end{remark}

\begin{remark}
We note that a preliminary version of the full classification of the
proper-biharmonic Legendre curves in Sasakian space forms was
obtained in ~\cite{Fetcu1}.
\end{remark}

In the following, we shall choose the unit $(2n+1)$-dimensional
sphere $\mathbb{S}^{2n+1}$ with its canonical and modified
Sasakian structures as a model for the complete, simply connected
Sasakian space form with constant $\varphi$-sectional curvature
$c>-3$, and we will find the explicit equations of biharmonic
Legendre curves obtained in the first three cases, viewed as
curves in $\mathbb{R}^{2n+2}$.

\begin{theorem}\label{curv2s2n+1,1}
Let $\gamma:I\to
(\mathbb{S}^{2n+1},\varphi_{0},\xi_{0},\eta_{0},g_{0})$, $n\geq
2$, be a proper-biharmonic Legendre curve parametrized by arc
length. Then the equation of $\gamma$ in the Euclidean space
$\mathbb{E}^{2n+2}=(\mathbb{R}^{2n+2},\langle,\rangle)$, is either
$$
\gamma(s)=\frac{1}{\sqrt{2}}\cos\Big(\sqrt{2}s\Big)e_{1}+\frac{1}{\sqrt{2}}\sin\Big(\sqrt{2}s\Big)e_{2}+\frac{1}{\sqrt{2}}e_{3}
$$
where $\{e_{i},\mathcal{I}e_{j}\}$ are constant unit vectors
orthogonal to each other, or
$$
\gamma(s)=\frac{1}{\sqrt{2}}\cos(As)e_{1}+\frac{1}{\sqrt{2}}\sin(As)e_{2}+\frac{1}{\sqrt{2}}\cos(Bs)e_{3}+\frac{1}{\sqrt{2}}\sin(Bs)e_{4},
$$
where
\begin{equation}\label{ecAB1}
A=\sqrt{1+\kappa_{1}},\ \ \ B=\sqrt{1-\kappa_{1}},\ \ \
\kappa_{1}\in(0,1),
\end{equation}
and $\{e_{i}\}$ are constant unit vectors orthogonal to each
other, with
$$
\langle e_{1},\mathcal{I}e_{3}\rangle=\langle
e_{1},\mathcal{I}e_{4}\rangle=\langle
e_{2},\mathcal{I}e_{3}\rangle=\langle
e_{2},\mathcal{I}e_{4}\rangle=0,\ \ A\langle
e_{1},\mathcal{I}e_{2}\rangle+B\langle
e_{3},\mathcal{I}e_{4}\rangle=0.
$$
\end{theorem}

\begin{proof} Let us denote by  $\dot\nabla$ and by $\widetilde{\nabla}$ the
Levi-Civita connections on $(\mathbb{S}^{2n+1},$ $g_{0})$ and
$(\mathbb{R}^{2n+2},\langle,\rangle)$, respectively.

First, assume that $\gamma$ is the biharmonic circle, that is
$\kappa_{1}=1$. From the Gauss and Frenet equations we get
$$
\widetilde{\nabla}_{T}T=\dot\nabla_{T}T-\langle
T,T\rangle\gamma=\kappa_{1}E_{2}-\gamma
$$
and
$$
\widetilde{\nabla}_{T}\widetilde{\nabla}_{T}T=(-\kappa_{1}^{2}-1)T=-2T,
$$
which implies
$$
\gamma^{'''}+2\gamma'=0.
$$
The general solution of the above equation is
$$
\gamma(s)=\cos\Big(\sqrt{2}s\Big)c_{1}+\sin\Big(\sqrt{2}s\Big)c_{2}+c_{3},
$$
where  $\{c_{i}\}$ are constant vectors in $\mathbb{E}^{2n+2}$.

\noindent Now, as $\gamma$ satisfies
$$
\langle\gamma,\gamma\rangle=1,\ \langle\gamma',\gamma'\rangle=1,\
\langle\gamma,\gamma'\rangle=0,\
\langle\gamma',\gamma''\rangle=0,\
\langle\gamma'',\gamma''\rangle=2,\
\langle\gamma,\gamma''\rangle=-1,
$$
and since in $s=0$ we have $\gamma=c_{1}+c_{3}$,
$\gamma'=\sqrt{2}c_{2}$, $\gamma''=-2c_{1}$, we obtain
$$
c_{11}+2c_{13}+c_{33}=1,\ c_{22}=\frac{1}{2},\ c_{12}+c_{23}=0,\
c_{12}=0,\ c_{11}=\frac{1}{2}, c_{11}+c_{13}=\frac{1}{2},
$$
where $c_{ij}=\langle c_{i},c_{j}\rangle$. Further, we get that
$\{c_{i}\}$ are orthogonal vectors in $\mathbb{E}^{2n+2}$ with
$\vert c_{1}\vert=\vert c_{2}\vert=\vert
c_{3}\vert=\frac{1}{\sqrt{2}}$.

\noindent Finally, using the fact that $\gamma$ is a Legendre
curve one obtains easily that $\langle
c_{i},\mathcal{I}c_{j}\rangle=0$ for any $i,j=1,2,3$.

Suppose now $\gamma$ is the biharmonic helix, that is
$\kappa_{1}^{2}+\kappa_{2}^{2}=1$, $\kappa_{1}\in(0,1)$. From the
Gauss and Frenet equations we get
$$
\widetilde{\nabla}_{T}T=\dot\nabla_{T}T-\langle
T,T\rangle\gamma=\kappa_{1}E_{2}-\gamma,
$$
$$
\widetilde{\nabla}_{T}\widetilde{\nabla}_{T}T=\kappa_{1}\widetilde{\nabla}_{T}E_{2}-T=\kappa_{1}\Big(-\kappa_{1}T+\kappa_{2}E_{3}\Big)-T=-\Big(\kappa_{1}^{2}+1\Big)T+\kappa_{1}\kappa_{2}E_{3}
$$
and
$$
\widetilde{\nabla}_{T}\widetilde{\nabla}_{T}\widetilde{\nabla}_{T}T=-\Big(\kappa_{1}^{2}+1\Big)\widetilde{\nabla}_{T}T+\kappa_{1}\kappa_{2}\widetilde{\nabla}_{T}E_{3}=
-\Big(\kappa_{1}^{2}+1\Big)\widetilde{\nabla}_{T}T-\kappa_{1}\kappa_{2}^{2}E_{2}=-2\gamma''-\kappa_{2}^{2}\gamma.
$$

\noindent Hence
$$
\gamma^{iv}+2\gamma''+\kappa_{2}^{2}\gamma=0,
$$
and its general solution is
$$
\gamma(s)=\cos(As)c_{1}+\sin(As)c_{2}+\cos(Bs)c_{3}+\sin(Bs)c_{4},
$$
where $A$, $B$ are given by (\ref{ecAB1}) and $\{c_{i}\}$ are
constant vectors in $\mathbb{E}^{2n+2}$.

\noindent As $\gamma$ satisfies
$$
\langle\gamma,\gamma\rangle=1,\ \langle\gamma',\gamma'\rangle=1,\
\langle\gamma,\gamma'\rangle=0,\
\langle\gamma',\gamma''\rangle=0,\
\langle\gamma'',\gamma''\rangle=1+\kappa_{1}^{2},
$$
$$
\langle\gamma,\gamma''\rangle=-1,\
\langle\gamma',\gamma'''\rangle=-(1+\kappa_{1}^{2}),\
\langle\gamma'',\gamma'''\rangle=0,\
\langle\gamma,\gamma'''\rangle=0,\
\langle\gamma''',\gamma'''\rangle=3\kappa_{1}^{2}+1,
$$
and since in $s=0$ we have $\gamma=c_{1}+c_{3}$,
$\gamma'=Ac_{2}+Bc_{4}$, $\gamma''=-A^{2}c_{1}-B^{2}c_{3}$,
$\gamma'''=-A^{3}c_{2}-B^{3}c_{4}$, we obtain
\begin{equation}\label{1.11}
c_{11}+2c_{13}+c_{33}=1
\end{equation}
\begin{equation}\label{1.21}
A^{2}c_{22}+2ABc_{24}+B^{2}c_{44}=1
\end{equation}
\begin{equation}\label{1.31}
Ac_{12}+Ac_{23}+Bc_{14}+Bc_{34}=0
\end{equation}
\begin{equation}\label{1.41}
A^{3}c_{12}+AB^{2}c_{23}+A^{2}Bc_{14}+B^{3}c_{34}=0
\end{equation}
\begin{equation}\label{1.51}
A^{4}c_{11}+2A^{2}B^{2}c_{13}+B^{4}c_{33}=1+\kappa_{1}^{2}
\end{equation}
\begin{equation}\label{1.61}
A^{2}c_{11}+(A^{2}+B^{2})c_{13}+B^{2}c_{33}=1
\end{equation}
\begin{equation}\label{1.71}
A^{4}c_{22}+(AB^{3}+A^{3}B)c_{24}+B^{4}c_{44}=1+\kappa_{1}^{2}
\end{equation}
\begin{equation}\label{1.81}
A^{5}c_{12}+A^{3}B^{2}c_{23}+A^{2}B^{3}c_{14}+B^{5}c_{34}=0
\end{equation}
\begin{equation}\label{1.91}
A^{3}c_{12}+A^{3}c_{23}+B^{3}c_{14}+B^{3}c_{34}=0
\end{equation}
\begin{equation}\label{1.101}
A^{6}c_{22}+2A^{3}B^{3}c_{24}+B^{6}c_{44}=3\kappa_{1}^{2}+1
\end{equation}
where $c_{ij}=\langle c_{i},c_{j}\rangle$. Since the determinant
of the system given by (\ref{1.31}), (\ref{1.41}), (\ref{1.81})
and (\ref{1.91}) is $-A^{2}B^{2}(A^{2}-B^{2})^{4}\neq 0$ it
follows that
$$
c_{12}=c_{23}=c_{14}=c_{34}=0.
$$

\noindent The equations (\ref{1.11}), (\ref{1.51}) and
(\ref{1.61}) give
$$
c_{11}=\frac{1}{2},\ \ c_{13}=0,\ \ c_{33}=\frac{1}{2},
$$
and, from (\ref{1.21}), (\ref{1.71}) and (\ref{1.101}) follows that
$$
c_{22}=\frac{1}{2},\ \ c_{24}=0,\ \ c_{44}=\frac{1}{2}.
$$

\noindent Therefore, we obtain that $\{c_{i}\}$ are orthogonal
vectors in $\mathbb{E}^{8}$ with $\vert c_{1}\vert=\vert
c_{2}\vert=\vert c_{3}\vert=\vert c_{4}\vert=\frac{1}{\sqrt{2}}$.

Finally, since $\gamma$ is a Legendre curve one obtains the
conclusion of the Theorem.
\end{proof}

\begin{theorem}\label{curv2s2n+1}
Let $\gamma:I\to (\mathbb{S}^{2n+1},\varphi,\xi,\eta,g)$, $n\geq 2$,
$a>0$, $a\neq 1$, be a proper-biharmonic Legendre curve parametrized
by arc length such that
$g(\nabla_{\gamma'}\gamma',\varphi\gamma')=0$. Then the equation of
$\gamma$ in the Euclidean space $\mathbb{E}^{2n+2}$, is either
$$
\gamma(s)=\frac{1}{\sqrt{2}}\cos\Big(\sqrt{\frac{2}{a}}s\Big)e_{1}+\frac{1}{\sqrt{2}}\sin\Big(\sqrt{\frac{2}{a}}s\Big)e_{2}+\frac{1}{\sqrt{2}}e_{3}
$$
for $n\geq 2$ or, for $n\geq 3$,
$$
\gamma(s)=\frac{1}{\sqrt{2}}\cos(As)e_{1}+\frac{1}{\sqrt{2}}\sin(As)e_{2}+\frac{1}{\sqrt{2}}\cos(Bs)e_{3}+\frac{1}{\sqrt{2}}\sin(Bs)e_{4},
$$
where
\begin{equation}\label{ecAB}
A=\sqrt{\frac{1+\kappa_{1}\sqrt{a}}{a}},\ \ \
B=\sqrt{\frac{1-\kappa_{1}\sqrt{a}}{a}},\ \ \
\kappa_{1}\in\Big(0,\frac{1}{a}\Big),
\end{equation}
and $\{e_{i},\mathcal{I}e_{j}\}$ are constant unit vectors
orthogonal to each other.
\end{theorem}

\begin{proof} Let us denote by  $\nabla$, $\dot\nabla$ and
by $\widetilde{\nabla}$ the Levi-Civita connections on
$(\mathbb{S}^{2n+1},g)$, $(\mathbb{S}^{2n+1},$ $g_{0})$ and
$(\mathbb{R}^{2n+2},\langle,\rangle)$, respectively.

First we consider the case when $\gamma$ is the biharmonic circle,
that is $\kappa_{1}^{2}=\frac{c+3}{4}$. Let $T=\gamma'$ be the unit
tangent vector field (with respect to the metric $g$) along
$\gamma$. Using the two Sasakian structures on $\mathbb{S}^{2n+1}$
we obtain $\dot\nabla_{T}T=\nabla_{T}T$ and
$\dot\nabla_{T}E_{2}=\nabla_{T}E_{2}$.

From the Gauss and Frenet equations we get
$$
\widetilde{\nabla}_{T}T=\dot\nabla_{T}T-\langle
T,T\rangle\gamma=\kappa_{1}E_{2}-\frac{1}{a}\gamma
$$
and
$$
\widetilde{\nabla}_{T}\widetilde{\nabla}_{T}T=(-\kappa_{1}^{2}-\frac{1}{a})T=-\frac{2}{a}T.
$$
Hence
$$
a\gamma^{'''}+2\gamma'=0,
$$
with the general solution
$$
\gamma(s)=\cos\Big(\sqrt{\frac{2}{a}}s\Big)c_{1}+\sin\Big(\sqrt{\frac{2}{a}}s\Big)c_{2}+c_{3},
$$
where  $\{c_{i}\}$ are constant vectors in $\mathbb{E}^{2n+2}$.

\noindent As $\gamma$ verifies the following equations,
$$
\langle\gamma,\gamma\rangle=1,\
\langle\gamma',\gamma'\rangle=\frac{1}{a},\
\langle\gamma,\gamma'\rangle=0,\
\langle\gamma',\gamma''\rangle=0,\
\langle\gamma'',\gamma''\rangle=\frac{2}{a^{2}},\
\langle\gamma,\gamma''\rangle=-\frac{1}{a},
$$
and in $s=0$ we have $\gamma=c_{1}+c_{3}$,
$\gamma'=\sqrt{\frac{2}{a}}c_{2}$, $\gamma''=-\frac{2}{a}c_{1}$,
one obtains
$$
c_{11}+2c_{13}+c_{33}=1,\ c_{22}=\frac{1}{2},\ c_{12}+c_{23}=0,\
c_{12}=0,\ c_{11}=\frac{1}{2}, c_{11}+c_{13}=\frac{1}{2},
$$
where $c_{ij}=\langle c_{i},c_{j}\rangle$. Consequently, we obtain
that $\{c_{i}\}$ are orthogonal vectors in $\mathbb{E}^{2n+2}$
with $\vert c_{1}\vert=\vert c_{2}\vert=\vert
c_{3}\vert=\frac{1}{\sqrt{2}}$.

\noindent Finally, using the facts that $\gamma$ is a Legendre
curve and $g(\nabla_{\gamma'}\gamma',\varphi\gamma')=0$ one
obtains easily that $\langle c_{i},\mathcal{I}c_{j}\rangle=0$ for
any $i,j=1,2,3$.

Now we assume that $\gamma$ is a biharmonic helix, that is
$\kappa_{1}^{2}+\kappa_{2}^{2}=\frac{c+3}{4}$,
$\kappa_{1}^{2}\in\Big(0,\frac{c+3}{4}\Big)$. First we obtain
$\dot\nabla_{T}T=\nabla_{T}T$,
$\dot\nabla_{T}E_{2}=\nabla_{T}E_{2}$ and
$\dot\nabla_{T}E_{3}=\nabla_{T}E_{3}$.

\noindent From the Gauss and Frenet equations we get
$$
\widetilde{\nabla}_{T}T=\dot\nabla_{T}T-\langle
T,T\rangle\gamma=\kappa_{1}E_{2}-\frac{1}{a}\gamma,
$$
$$
\widetilde{\nabla}_{T}\widetilde{\nabla}_{T}T=\kappa_{1}\widetilde{\nabla}_{T}E_{2}-\frac{1}{a}T=
\kappa_{1}\Big(-\kappa_{1}T+\kappa_{2}E_{3}\Big)-\frac{1}{a}T=-\Big(\kappa_{1}^{2}+\frac{1}{a}\Big)T+\kappa_{1}\kappa_{2}E_{3},
$$
and
$$
\begin{array}{ll}
\widetilde{\nabla}_{T}\widetilde{\nabla}_{T}\widetilde{\nabla}_{T}T&=-\Big(\kappa_{1}^{2}+\frac{1}{a}\Big)\widetilde{\nabla}_{T}T+\kappa_{1}\kappa_{2}\widetilde{\nabla}_{T}E_{3}=
-\Big(\kappa_{1}^{2}+\frac{1}{a}\Big)\widetilde{\nabla}_{T}T-\kappa_{1}\kappa_{2}^{2}E_{2}\\
\\&=-\frac{2}{a}\gamma''-\frac{1}{a}\kappa_{2}^{2}\gamma.
\end{array}
$$

\noindent Therefore
$$
a\gamma^{iv}+2\gamma''+\kappa_{2}^{2}\gamma=0,
$$
and its general solution is
$$
\gamma(s)=\cos(As)c_{1}+\sin(As)c_{2}+\cos(Bs)c_{3}+\sin(Bs)c_{4},
$$
where $A$, $B$ are given by (\ref{ecAB}) and $\{c_{i}\}$ are
constant vectors in $\mathbb{E}^{2n+2}$.

\noindent The curve $\gamma$ satisfies
$$
\langle\gamma,\gamma\rangle=1,\
\langle\gamma',\gamma'\rangle=\frac{1}{a},\
\langle\gamma,\gamma'\rangle=0,\
\langle\gamma',\gamma''\rangle=0,\
\langle\gamma'',\gamma''\rangle=\frac{1+a\kappa_{1}^{2}}{a^{2}},
$$
$$
\langle\gamma,\gamma''\rangle=-\frac{1}{a},\
\langle\gamma',\gamma'''\rangle=-\frac{1+a\kappa_{1}^{2}}{a^{2}},\
\langle\gamma'',\gamma'''\rangle=0,\
\langle\gamma,\gamma'''\rangle=0,\
$$
$$
\langle\gamma''',\gamma'''\rangle=\frac{3a\kappa_{1}^{2}+1}{a^{3}},
$$
and in $s=0$ we have
$$
\gamma=c_{1}+c_{3},\ \ \gamma'=Ac_{2}+Bc_{4},\ \
\gamma''=-A^{2}c_{1}-B^{2}c_{3},\ \
\gamma'''=-A^{3}c_{2}-B^{3}c_{4}.
$$
Then, it follows
\begin{equation}\label{1.1}
c_{11}+2c_{13}+c_{33}=1
\end{equation}
\begin{equation}\label{1.2}
A^{2}c_{22}+2ABc_{24}+B^{2}c_{44}=\frac{1}{a}
\end{equation}
\begin{equation}\label{1.3}
Ac_{12}+Ac_{23}+Bc_{14}+Bc_{34}=0
\end{equation}
\begin{equation}\label{1.4}
A^{3}c_{12}+AB^{2}c_{23}+A^{2}Bc_{14}+B^{3}c_{34}=0
\end{equation}
\begin{equation}\label{1.5}
A^{4}c_{11}+2A^{2}B^{2}c_{13}+B^{4}c_{33}=\frac{1+a\kappa_{1}^{2}}{a^{2}}
\end{equation}
\begin{equation}\label{1.6}
A^{2}c_{11}+(A^{2}+B^{2})c_{13}+B^{2}c_{33}=\frac{1}{a}
\end{equation}
\begin{equation}\label{1.7}
A^{4}c_{22}+(AB^{3}+A^{3}B)c_{24}+B^{4}c_{44}=\frac{1+a\kappa_{1}^{2}}{a^{2}}
\end{equation}
\begin{equation}\label{1.8}
A^{5}c_{12}+A^{3}B^{2}c_{23}+A^{2}B^{3}c_{14}+B^{5}c_{34}=0
\end{equation}
\begin{equation}\label{1.9}
A^{3}c_{12}+A^{3}c_{23}+B^{3}c_{14}+B^{3}c_{34}=0
\end{equation}
\begin{equation}\label{1.10}
A^{6}c_{22}+2A^{3}B^{3}c_{24}+B^{6}c_{44}=\frac{3a\kappa_{1}^{2}+1}{a^{3}}
\end{equation}
where $c_{ij}=\langle c_{i},c_{j}\rangle$.

\noindent The solution of the system given by (\ref{1.3}),
(\ref{1.4}), (\ref{1.8}) and (\ref{1.9}) is
$$
c_{12}=c_{23}=c_{14}=c_{34}=0.
$$

\noindent From equations (\ref{1.1}), (\ref{1.5}) and (\ref{1.6})
we get
$$
c_{11}=\frac{1}{2},\ \ c_{13}=0,\ \ c_{33}=\frac{1}{2},
$$
and, from (\ref{1.2}), (\ref{1.7}), (\ref{1.10}),
$$
c_{22}=\frac{1}{2},\ \ c_{24}=0,\ \ c_{44}=\frac{1}{2}.
$$

\noindent We obtain that $\{c_{i}\}$ are orthogonal vectors in
$\mathbb{E}^{8}$ with $\vert c_{1}\vert=\vert c_{2}\vert=\vert
c_{3}\vert=\vert c_{4}\vert=\frac{1}{\sqrt{2}}$.

\noindent Finally, since $\gamma$ is a Legendre curve and
$g(\nabla_{\gamma'}\gamma',\varphi\gamma')=0$, one obtains the
conclusion.

\end{proof}

Just like for $\mathbb{S}^{3}$ (see \cite{FetcuOniciuc}) one
obtains for $\mathbb{S}^{2n+1}$ endowed with the modified Sasakian
structure, with $0<a<1$ (that means $c>1$), the following result.

\begin{theorem}\label{curv1s2n+1}
Let $\gamma:I\to (\mathbb{S}^{2n+1},\varphi,\xi,\eta,g)$ be a
biharmonic Legendre curve parametrized by its arc length such that
$\nabla_{\gamma'}\gamma'=\sqrt{c-1}\varphi\gamma'$. Then the
equation of $\gamma$ in the Euclidean space $\mathbb{E}^{2n+2}$ is
$$
\begin{array}{rlc}
\gamma(s)&=&\sqrt{\frac{B}{A+B}}\cos(As)e_{1}-\sqrt{\frac{B}{A+B}}\sin(As)\mathcal{I}e_{1}\\
\\&&+\sqrt{\frac{A}{A+B}}\cos(Bs)e_{3}+\sqrt{\frac{A}{A+B}}\sin(Bs)\mathcal{I}{e_{3}}\\ \\&=&\sqrt{\frac{B}{A+B}}\exp(-\mathrm{i}As)e_{1}+\sqrt{\frac{A}{A+B}}\exp(\mathrm{i}Bs)e_{3},
\end{array}
$$
where $\{e_{i}\}$ are constant unit orthogonal vectors in
$\mathbb{E}^{2n+2}$ and
\begin{equation}\label{4}
A=\sqrt{\frac{3-2a-2\sqrt{(a-1)(a-2)}}{a}},\ \
B=\sqrt{\frac{3-2a+2\sqrt{(a-1)(a-2)}}{a}}.
\end{equation}
\end{theorem}

\begin{remark}The ODE satisfied by proper-biharmonic Legendre curves in the
$(2n+1)$-sphere, in the fourth case, may be also obtained but the
computations are rather complicated.
\end{remark}

\section{Biharmonic submanifolds in Sasakian space forms}
\setcounter{equation}{0}

A method to obtain biharmonic submanifolds in a Sasakian space
form is provided by the following Theorem.

\begin{theorem} \label{teorema1}
Let $(N^{2n+1},\varphi,\xi,\eta,g)$ be a strictly regular Sasakian
space form with constant $\varphi$-sectional curvature $c$ and let
${\bf i}:M\to N$ be an $r$-dimensional integral submanifold of $N$.
Consider
$$
F:\widetilde{M}=I\times M\to N,\ \ \
F(t,p)=\phi_{t}(p)=\phi_{p}(t),
$$
where $I=\mathbb{S}^{1}$ or $I=\mathbb{R}$ and
$\{\phi_{t}\}_{t\in\mathbb{R}}$ is the flow of the vector field
$\xi$. Then $F:(\widetilde{M},\widetilde{g}=dt^{2}+{\bf
i}^{\ast}g)\to N$ is a Riemannian immersion and it is
proper-biharmonic if and only if $M$ is a proper-biharmonic
submanifold of $N$.
\end{theorem}

\begin{proof}From the definition of the flow of $\xi$ we have
$$
dF(t,p)\Big(\frac{\partial}{\partial
t}\Big)=\frac{d}{ds}\vert_{s=t}\{\phi_{p}(s)\}=\dot\phi_{p}(t)=\xi(\phi_{p}(t))=\xi(F(t,p)),
$$
i.e. $\frac{\partial}{\partial t}$ is $F$-correlated to $\xi$ and
$$
\Big\vert dF(t,p)\Big(\frac{\partial}{\partial
t}\Big)\Big\vert=\vert \xi(F(t,p))\vert=1=\Big\vert
\frac{\partial}{\partial t}\Big\vert.
$$

\noindent The vector $X_{p}\in T_{p}M$ can be identified to
$(0,X_{p})\in T_{(t',p)}(I\times M)$, for any $t'\in I$, and we
have
$$
dF_{(t,p)}(X_{p})=(dF)_{(t,p)}(\dot\gamma(0))=\frac{d}{ds}\vert_{s=0}\{\phi_{t}(\gamma(s))\}=(d\phi_{t})_{p}(X_{p}).
$$
Since $\phi_{t}$ is an isometry, we have $\vert
dF_{(t,p)}(X_{p})\vert=\vert (d\phi_{t})_{p}(X_{p})\vert=\vert
X_{p}\vert$.

\noindent Moreover,
$$
\begin{array}{ll}
g\Big(dF_{(t,p)}\Big(\frac{\partial}{\partial
t}\Big),dF_{(t,p)}(X_{p})\Big)&=g(\xi(\phi_{p}(t)),(d\phi_{t})_{p}(X_{p}))\\
\\&=g((d\phi_{t})_{p}(\xi_{p}),(d\phi_{t})_{p}(X_{p}))=g(\xi_{p},X_{p})\\
\\&=0,
\end{array}
$$
so $F:(I\times M,dt^{2}+{\bf i}^{\ast}g)\to N$ is a Riemannian
immersion.

Let  $F^{-1}(TN)$ be the pull-back bundle over $\widetilde{M}$ and
$\nabla^{F}$ the pull-back connection determined by the
Levi-Civita connection on $N$. We shall prove that
$$
\tau(F)_{(t,p)}=(d\phi_{t})_{p}(\tau(\textbf{i}))\ \
\textnormal{and}\ \
\tau_{2}(F)_{(t,p)}=(d\phi_{t})_{p}(\tau_{2}(\textbf{i})),
$$
so, from the point of view of harmonicity and biharmonicity,
$\widetilde{M}$ and $M$ have the same behaviour.

We start with two remarks. Let $\sigma\in F^{-1}(TN)$ defined by
$\sigma_{(t,p)}=(d\phi_{t})_{p}(Z_{p})$, where $Z$ is a vector
field along $M$, that is $Z_{p}\in T_{p}N$, $\forall p\in M$. We
have
\begin{equation}\label{ecfi1}
(\nabla^{F}_{X}\sigma)_{(t,p)}=(d\phi_{t})_{p}(\nabla^{N}_{X}Z),\
\ \forall X\in C(TM).
\end{equation}

Then, if $\sigma\in F^{-1}(TN)$, it follows that $\varphi\sigma\in
F^{-1}(TN)$,
$(\varphi\sigma)_{(t,p)}=\varphi_{\phi_{p}(t)}(\sigma_{(t,p)})$,
and
\begin{equation}\label{ecfi2}
\nabla^{F}_{\frac{\partial}{\partial
t}}\varphi\sigma=\varphi\nabla^{F}_{\frac{\partial}{\partial
t}}\sigma.
\end{equation}

\noindent Now, we consider $\{X_{1},...,X_{r}\}$ a local
orthonormal frame field on $U$, where $U$ is an open subset of
$M$. The tension field of $F$ is given by
$$
\tau(F)=\nabla^{F}_{\frac{\partial}{\partial
t}}dF\Big(\frac{\partial}{\partial
t}\Big)-dF\Big(\nabla^{\widetilde{M}}_{\frac{\partial}{\partial
t}}\frac{\partial}{\partial
t}\Big)+\sum_{a=1}^{r}\{\nabla^{F}_{X_{a}}dF(X_{a})-dF(\nabla^{\widetilde{M}}_{X_{a}}X_{a})\}.
$$
As
$$
\nabla^{F}_{\frac{\partial}{\partial
t}}dF\Big(\frac{\partial}{\partial t}\Big)=\nabla^{N}_{\xi}\xi=0,\
\  \nabla^{\widetilde{M}}_{\frac{\partial}{\partial
t}}\frac{\partial}{\partial
t}=\nabla^{I}_{\frac{\partial}{\partial
t}}\frac{\partial}{\partial t}=0,
$$
$$
(\nabla^{F}_{X_{a}}dF(X_{a}))_{(t,p)}=(d\phi_{t})_{p}(\nabla^{N}_{X_{a}}X_{a}),\
\
dF_{(t,p)}(\nabla^{\widetilde{M}}_{X_{a}}X_{a})=(d\phi_{t})_{p}(\nabla^{M}_{X_{a}}X_{a})
$$
we obtain
$$
\tau(F)_{(t,p)}=(d\phi_{t})_{p}(\tau(\textbf{i})).
$$

\noindent The next step is to prove that
$\nabla^{F}_{\frac{\partial}{\partial
t}}\tau(F)=-\varphi(\tau(F))$. Since $[\frac{\partial}{\partial
t},X_{a}]=0$, $a=1,...,r$,
$$
\nabla^{F}_{\frac{\partial}{\partial
t}}dF(X_{a})=\nabla^{F}_{X_{a}}dF\Big(\frac{\partial}{\partial
t}\Big).
$$
But
$$
\begin{array}{ll}
\Big(\nabla^{F}_{X_{a}}dF\Big(\frac{\partial}{\partial
t}\Big)\Big)_{(t,p)}&=\nabla^{N}_{dF_{(t,p)}X_{a}}\xi=\nabla^{N}_{(d\phi_{t})_{p}X_{a}}\xi=-\varphi((d\phi_{t})_{p}(X_{a}))\\
\\&=-(d\phi_{t})_{p}(\varphi X_{a}),\end{array}
$$
so
\begin{equation}\label{ecfi3}
\Big(\nabla^{F}_{\frac{\partial}{\partial
t}}dF(X_{a})\Big)_{(t,p)}=-(d\phi_{t})_{p}(\varphi X_{a}).
\end{equation}

\noindent We note that
$$
R^{F}\Big(\frac{\partial}{\partial
t},X_{a}\Big)dF(X_{a})=\nabla^{F}_{\frac{\partial}{\partial
t}}\nabla^{F}_{X_{a}}dF(X_{a})-\nabla^{F}_{X_{a}}\nabla^{F}_{\frac{\partial}{\partial
t}}dF(X_{a})
$$
and, on the other hand
$$
\Big(R^{F}\Big(\frac{\partial}{\partial
t},X_{a}\Big)dF(X_{a})\Big)_{(t,p)}=R^{N}_{\phi_{t}(p)}(\xi,(d\phi_{t})_{p}(X_{a}))(d\phi_{t})_{p}(X_{a})=\xi.
$$

\noindent Therefore
\begin{equation}\label{ecfi6}
\nabla^{F}_{\frac{\partial}{\partial
t}}\nabla^{F}_{X_{a}}dF(X_{a})-\nabla^{F}_{X_{a}}\nabla^{F}_{\frac{\partial}{\partial
t}}dF(X_{a})=\xi
\end{equation}

\noindent From (\ref{ecfi1}) and (\ref{ecfi3}) we get
\begin{equation}\label{ecfi7}\begin{array}{ll}
\Big(\nabla^{F}_{X_{a}}\nabla^{F}_{\frac{\partial}{\partial
t}}dF(X_{a})\Big)_{(t,p)}&=-(d\phi_{t})_{p}(\nabla^{N}_{X_{a}}\varphi
X_{a})\\
\\&=-(d\phi_{t})_{p}(\xi+\varphi\nabla^{N}_{X_{a}}X_{a})\end{array}
\end{equation}
and from (\ref{ecfi3})
\begin{equation}\label{ecfi8}\begin{array}{ll}
\Big(\nabla^{F}_{\frac{\partial}{\partial
t}}dF(\nabla^{\widetilde{M}}_{X_{a}}X_{a})\Big)_{(t,p)}&=\Big(\nabla^{F}_{\frac{\partial}{\partial
t}}dF(\nabla^{M}_{X_{a}}X_{a})\Big)_{(t,p)}\\
\\&=-(d\phi_{t})_{p}(\varphi\nabla^{M}_{X_{a}}X_{a}).\end{array}
\end{equation}

\noindent Replacing (\ref{ecfi7}) and (\ref{ecfi8}) in
(\ref{ecfi6}), we obtain
$$
\begin{array}{lll}
\xi&=&\nabla^{F}_{\frac{\partial}{\partial
t}}\nabla^{F}_{X_{a}}dF(X_{a})-\nabla^{F}_{\frac{\partial}{\partial
t}}dF(\nabla^{\widetilde{M}}_{X_{a}}X_{a})+\nabla^{F}_{\frac{\partial}{\partial
t}}dF(\nabla^{\widetilde{M}}_{X_{a}}X_{a})-\nabla^{F}_{X_{a}}\nabla^{F}_{\frac{\partial}{\partial
t}}dF(X_{a})\\ \\&=&\nabla^{F}_{\frac{\partial}{\partial t}}\nabla
dF(X_{a},X_{a})-(d\phi_{t})_{p}(\varphi\nabla^{M}_{X_{a}}X_{a})+(d\phi_{t})_{p}(\xi+\varphi\nabla^{N}_{X_{a}}X_{a}),
\end{array}
$$
so
$$
\Big(\nabla^{F}_{\frac{\partial}{\partial
t}}\nabla^{F}_{X_{a}}dF(X_{a})\Big)_{(t,p)}=-\varphi(d\phi_{t})_{p}(\nabla
d\textbf{i}(X_{a},X_{a})).
$$
As $\nabla dF(X_{a},X_{a})=0$, we obtain
$$
\nabla^{F}_{\frac{\partial}{\partial t}}\tau(F)=-\varphi(\tau(F)).
$$
From (\ref{ecfi2}) we have
\begin{equation}\label{ecfi9}
\begin{array}{ll}
\nabla^{F}_{\frac{\partial}{\partial
t}}\nabla^{F}_{\frac{\partial}{\partial
t}}\tau(F)&=-\nabla^{F}_{\frac{\partial}{\partial
t}}\varphi(\tau(F))=-\varphi\nabla^{F}_{\frac{\partial}{\partial
t}}\tau(F)=\varphi^{2}\tau(F)\\ \\&=-\tau(F),\end{array}
\end{equation}
and from (\ref{ecfi1})
\begin{equation}\label{ecfi10}
(\nabla^{F}_{X_{a}}\nabla^{F}_{X_{a}}\tau(F))_{(t,p)}=(d\phi_{t})_{p}(\nabla^{N}_{X_{a}}\nabla^{N}_{X_{a}}\tau(\textbf{i})),
\end{equation}
\begin{equation}\label{ecfi11}
\Big(\nabla^{F}_{\nabla^{\widetilde{M}}_{X_{a}}X_{a}}\tau(F)\Big)_{(t,p)}=(d\phi_{t})_{p}\Big(\nabla^{N}_{\nabla^{M}_{X_{a}}X_{a}}\tau(\textbf{i})\Big).
\end{equation}

\noindent From (\ref{ecfi9}), (\ref{ecfi10}) and (\ref{ecfi11}) we
obtain
\begin{equation}\label{ec-delta}
-(\Delta^{F}\tau(F))_{(t,p)}=-\tau(F)_{(t,p)}-(d\phi_{t})_{p}(\Delta^{\textbf{i}}\tau(\textbf{i})).
\end{equation}

\noindent After a straightforward computation we get
\begin{equation}\label{ec-curbura}
\trace R^{F}(dF,\tau(F))dF=-\tau(F)+(d\phi_{t})_{p}(\trace
R^{N}_{p}(d\textbf{i},\tau(\textbf{i}))d\textbf{i}).
\end{equation}

\noindent Finally, from (\ref{ec-delta}) and (\ref{ec-curbura}) we
obtain
$$
\tau_{2}(F)_{(t,p)}=(d\phi_{t})_{p}(\tau_{2}(\textbf{i})).
$$

\end{proof}

\begin{remark}

The previous result was expected because of the following simple
remark. Assume that $(N^{2n+1},\varphi,\xi,\eta,g)$ is a compact
strictly regular Sasakian space form with constant
$\varphi$-sectional curvature $c$ and let $G:M\to N$ be an arbitrary
smooth map from a compact Riemannian manifold $M$. If $F$ is
biharmonic, then the map $G$ is biharmonic, where
$F:\widetilde{M}=\s^1\times M\to N$, $F(t,p)=\phi_t(G(p))$.

\noindent Indeed, an arbitrary variation $\{G_s\}_{s}$ of $G$
induces a variation $\{F_s\}_{s}$ of $F$. We have that
$\tau_{(p,t)}(F_s)=(d\phi_t)_{G_s(p)}(\tau_p(G_s))$ and, from the
biharmonicity of $F$ and the Fubini Theorem, we get
\begin{eqnarray*}
0&=&\frac{d}{ds}\vert_{s=0}\{E_2(F_s)\}=\frac{1}{2}\frac{d}{ds}\vert_{s=0}\int_{\widetilde{M}}\vert
\tau(F_s)\vert^2 \
v_{\widetilde{g}}=\frac{\pi}{2}\frac{d}{ds}\vert_{s=0}\int_M\vert\tau(G_s)\vert^2
\ v_g \\
&=&\pi\frac{d}{ds}\vert_{s=0}\{E_2(G_s)\}.
\end{eqnarray*}

\end{remark}

Next, consider the unit $(2n+1)$-dimensional sphere
$\mathbb{S}^{2n+1}$ endowed with its canonical or modified Sasakian
structure. The flow of $\xi$ is
$\phi_{t}(z)=\exp(-\mathrm{i}\frac{t}{a})z$, where
$$
z=(z^{1},...,z^{n+1})=(x^{1},...,x^{n+1},y^{1},...,y^{n+1}).
$$

\noindent From the above expression of the flow and from Theorems
\ref{curv1s2n+1}, \ref{curv2s2n+1}, \ref{curv2s2n+1,1} and
\ref{teorema1} we obtain explicit examples of proper-biharmonic
submanifolds in $(\mathbb{S}^{2n+1},\varphi,\xi,\eta,g)$, $a>0$,
of constant mean curvature. In particular, we reobtain a result in
\cite{Sasahara2} and, for $n=1$, the result in
\cite{FetcuOniciuc}.

\begin{proposition}[\cite{Sasahara2}] Let $F:\widetilde{M}\to \mathbb{S}^{5}\subset\mathbb{R}^{6}$ be a
proper-biharmonic anti-invariant immersion. Then the position
vector of $\widetilde{M}$ in $\mathbb{R}^{6}$ is
$$
F(t,u,v)=\frac{\exp(-\mathrm{i}t)}{\sqrt{2}}(\exp(\mathrm{i}u),\mathrm{i}\exp(-\mathrm{i}u)\sin(\sqrt{2}v),\mathrm{i}\exp(-\mathrm{i}u)\cos(\sqrt{2}v)).
$$
\end{proposition}

\begin{proof} It was proved in \cite{Sasahara1} that the
proper-biharmonic integral surface of $(\mathbb{S}^{5},$
$\varphi_{0},\xi_{0},\eta_{0},$ $g_{0})$ is given by
$$
f(u,v)=\frac{1}{\sqrt{2}}(\exp(\mathrm{i}u),\mathrm{i}\exp(-\mathrm{i}u)\sin\sqrt{2}v,\mathrm{i}\exp(-\mathrm{i}u)\cos\sqrt{2}v).
$$
Now, composing with the flow of $\xi_{0}$ we reobtain the result
in \cite{Sasahara2}.
\end{proof}

\begin{proposition}[\cite{FetcuOniciuc}]\label{surf1s2n+1} Let $M$ be a surface
in $(\mathbb{S}^{2n+1},\varphi,\xi,\eta,g)$ with $a\in\mathbb(0,1)$,
with the position vector in the Euclidean space $\mathbb{E}^{2n+2}$
given by
$$
F(t,s)=\exp\Big(-\mathrm{i}\frac{t}{a}\Big)\Big(\sqrt{\frac{B}{A+B}}\exp(-\mathrm{i}As)e_{1}+\sqrt{\frac{A}{A+B}}\exp(\mathrm{i}Bs)e_{3}\Big),
$$
where $\{e_{1},e_{3}\}$ is an orthonormal system of constant vectors
in the Euclidean space $(\mathbb{R}^{2n+2},\langle,\rangle)$, with
$e_{3}$ orthogonal to $\mathcal{I}e_{1}$ and $A$, $B$ are given by
(\ref{4}).

\noindent Then $M$ is a proper-biharmonic surface in
$\mathbb{S}^{2n+1}(c)$.
\end{proposition}

\section{Proper-biharmonic submanifolds of
$(\mathbb{S}^{7},g_{0})$} \setcounter{equation}{0}

We consider the Euclidean space $\mathbb{E}^{8}$ with three complex
structures,
$$
\mathcal{I}=\left(\begin{array}{cc}0&-I_{4}\\
I_{4}&0\end{array}\right),\ \
\ \mathcal{J}=\left(\begin{array}{cccc}0&0&0&I_{2}\\ 0&0&-I_{2}&0\\
0&I_{2}&0&0\\ -I_{2}&0&0&0\end{array}\right),\ \ \
\mathcal{K}=-\mathcal{I}\mathcal{J},
$$
where $I_{n}$ denotes the $n\times n$ identity matrix. We define
three vector fields on $\mathbb{S}^{7}$ by
$$
\xi_{1}=-\mathcal{I}z,\ \ \  \ \xi_{2}=-\mathcal{J}z, \ \ \ \
\xi_{3}=-\mathcal{K}z,\ \ \ z\in\mathbb{S}^{7},
$$
and consider their dual 1-forms $\eta_{1}=\eta_{0}$, $\eta_{2}$,
$\eta_{3}$. Let $\varphi_{a}$ defined by
$$
\varphi_{1}=\varphi_{0}=s\circ\mathcal{I},\ \
\varphi_{2}=s\circ\mathcal{J},\ \ \varphi_{3}=s\circ\mathcal{K}.
$$
Then $(\varphi_{a},\xi_{a},\eta_{a},g_{0})$, $a=1,2,3$, determine
a Sasakian 3-structure on $\mathbb{S}^{7}$.

Now, we shall indicate a method to construct proper-biharmonic
submanifolds in $(\mathbb{S}^{7},g_{0})$. We consider
$\gamma=\gamma(s)$ a proper-biharmonic curve in
$(\mathbb{S}^{7},g_{0})$, parametrized by arc-length, which is a
Legendre curve for two of the three contact structures (it was
proved in \cite{Fetcu} that there is no proper-biharmonic curve
which is Legendre with respect to all three contact structures on
$\mathbb{S}^{7}$). For example, assume that $\gamma$ is a Legendre
curve for $\eta_{1}$ and $\eta_{2}$. Composing with the flow of
$\xi_{1}$ (or $\xi_{2}$) we obtain a biharmonic surface which is
Legendre with respect to $\eta_{2}$ (or $\eta_{1}$). Then,
composing with the flow of $\xi_{2}$ (or $\xi_{1}$) we get a
biharmonic 3-dimensional submanifold of $(\mathbb{S}^{7},g_{0})$.

\noindent Using this method, from Theorems \ref{curv2s2n+1,1} and
\ref{teorema1}, we obtain 4 classes of proper-biharmonic surfaces
in $(\mathbb{S}^{7},g_{0})$ and 4 classes of proper-biharmonic
3-dimensional submanifolds of $(\mathbb{S}^{7},g_{0})$, all of
constant mean curvature.

\noindent For example, from Theorems \ref{curv2s2n+1,1} and
\ref{teorema1}, composing first with the flow of $\xi_{1}$ and
then with that of $\xi_{2}$, we get the explicit parametric
equations of proper-biharmonic 3-dimensional submanifolds of
$(\mathbb{S}^{7},g_{0})$.

\begin{proposition}\label{dim3s7} Let $M$ be a $3$-dimensional
submanifold in $\mathbb{S}^{7}$ such that its position vector
field in $\mathbb{E}^{8}$ is either
$$
\begin{array}{lll}
\mathbf{x}_{1}&=&\mathbf{x}_{1}(u,t,s)\\ &=&\frac{1}{\sqrt{2}}\Big(\cos(u)\cos(\sqrt{2}s)\cos(t)e_{1}+\cos(u)\sin(\sqrt{2}s)\cos(t)e_{2}\\
\\&&+\cos(u)\cos(t)e_{3}
-\cos(u)\cos(\sqrt{2}s)\sin(t)\mathcal{I}e_{1}\\
\\&&-\cos(u)\sin(\sqrt{2}s)\sin(t)\mathcal{I}e_{2}-\cos(u)\sin(t)\mathcal{I}e_{3}\\
\\&&-\sin(u)\cos(\sqrt{2}s)\cos(t)\mathcal{J}e_{1}-\sin(u)\sin(\sqrt{2}s)\cos(t)\mathcal{J}e_{2}\\
\\&&-\sin(u)\cos(t)\mathcal{J}e_{3}
-\sin(u)\cos(\sqrt{2}s)\sin(t)\mathcal{K}e_{1}\\
\\&&-\sin(u)\sin(\sqrt{2}s)\sin(t)\mathcal{K}e_{2}-\sin(u)\sin(t)\mathcal{K}e_{3}\Big),
\end{array}
$$
where $\{e_{i},\mathcal{I}e_{j}\}$, $\{e_{i},\mathcal{J}e_{j}\}$
are systems of constant orthonormal vectors in $\mathbb{E}^{8}$,
or
$$
\begin{array}{lll}
\mathbf{x}_{2}&=&\mathbf{x}_{2}(u,t,s)
\\&=&\frac{1}{\sqrt{2}}\Big(\cos(u)\cos(As)\cos(t)e_{1}+\cos(u)\sin(As)\cos(t)e_{2}\\
\\&&+\cos(u)\cos(Bs)\cos(t)e_{3}
+\cos(u)\sin(Bs)\cos(t)e_{4}\\
\\&&-\cos(u)\cos(As)\sin(t)\mathcal{I}e_{1}-\cos(u)\sin(As)\sin(t)\mathcal{I}e_{2}\\
\\&&-\cos(u)\cos(Bs)\sin(t)\mathcal{I}e_{3}-\cos(u)\sin(Bs)\sin(t)\mathcal{I}e_{4}\\
\\&&-\sin(u)\cos(As)\cos(t)\mathcal{J}e_{1}-\sin(u)\sin(As)\cos(t)\mathcal{J}e_{2}\\
\\&&-\sin(u)\cos(Bs)\cos(t)\mathcal{J}e_{3}
-\sin(u)\sin(Bs)\cos(t)\mathcal{J}e_{4}\\
\\
&&-\sin(u)\cos(As)\sin(t)\mathcal{K}e_{1}-\sin(u)\sin(As)\sin(t)\mathcal{K}e_{2}\\
\\&&-\sin(u)\cos(Bs)\sin(t)\mathcal{K}e_{3}-\sin(u)\sin(Bs)\sin(t)\mathcal{K}e_{4}\Big),
\end{array}
$$
where
$$
A=\sqrt{1+\kappa_{1}},\ \ \ B=\sqrt{1-\kappa_{1}},\ \ \
\kappa_{1}=\cst\in(0,1),
$$
and $\{e_{i}\}$ are constant orthonormal vectors in
$\mathbb{E}^{8}$ such that
$$
\langle e_{1},\mathcal{I}e_{3}\rangle=\langle
e_{1},\mathcal{I}e_{4}\rangle=\langle
e_{2},\mathcal{I}e_{3}\rangle=\langle
e_{2},\mathcal{I}e_{4}\rangle=0,
$$
$$
\langle e_{1},\mathcal{J}e_{3}\rangle=\langle
e_{1},\mathcal{J}e_{4}\rangle=\langle
e_{2},\mathcal{J}e_{3}\rangle=\langle
e_{2},\mathcal{J}e_{4}\rangle=0,
$$
and
$$
A\langle e_{1},\mathcal{I}e_{2}\rangle+B\langle
e_{3},\mathcal{I}e_{4}\rangle=A\langle
e_{1},\mathcal{J}e_{2}\rangle+B\langle
e_{3},\mathcal{J}e_{4}\rangle=0.
$$

\noindent Then $M$ is a proper-biharmonic submanifold of
$(\mathbb{S}^{7},g_{0})$.
\end{proposition}

\begin{proof} As the flows of $\xi_{1}$ and $\xi_{2}$ are given by
$$
\phi_{t}^{1}(z)=(\cos t)z-(\sin t)\mathcal{I}z,\ \ \ \
\phi_{t}^{2}(z)=(\cos t)z-(\sin t)\mathcal{J}z,
$$
the Proposition follows by a straightforward computation.

\end{proof}

\begin{remark} Note that there exist vectors $\{e_{i}\}$ which
satisfy the hypotheses of the previous Proposition. For example the
first three or four vectors, respectively, from the canonical basis
of $\mathbb{E}^{8}$ satisfy the required properties.
\end{remark}

\end{document}